\documentclass[10pt]{amsart}
	\usepackage[all]{xy}
	\usepackage{amssymb}
	\usepackage{fancybox}

	\newcommand{\T}{{\mathcal{T}}}
	\newcommand{\A}{{ \rm Aut }}

	\renewcommand{\O}{{\mathcal{O}}}

	\newcommand{\lf}{\left\lfloor}
	\newcommand{\rf}{\right\rfloor}

	\renewcommand{\mod}{{\;\rm mod}}

	\newcommand{\be}{\begin{equation}}
	\newcommand{\ee}{\end{equation}}

	\newtheorem{pro}{Proposition}[section]
	\newtheorem{lemma}[pro]{Lemma}

	\newtheorem{theorem}[pro]{Theorem}
	
	\newtheorem{definition}[pro]{Definition}

\begin{document}
	\bibliographystyle{amsplain}
	
\title{
Polydifferentials and the deformation functor of curves with 
	automorphisms II}

	\author{A. Kontogeorgis}

	\begin{abstract}
We apply the known results on the Galois module structure of the sheaf of polydifferentials 
in order to study the dimension of the tangent space of the deformation functor of curves 
with automorphisms. We are able to find the dimension  for the case of  weakly ramified covers 
and for the case of the action of a cyclic group of order $p^v$.
        \end{abstract}

	\address{
	Department of Mathematics, University of the \AE gean, 83200 Karlovassi, Samos,
	Greece\\ { \texttt{\upshape http://eloris.samos.aegean.gr}}
	}
\email{kontogar@aegean.gr}
	\date{\today}

	\maketitle

	\section{Introduction}

	Let $X$ be a  non-singular complete curve of genus $g\geq 2$ 
	defined over an algebraic closed field $k$ 
	of positive characteristic $p>3$, and let $G$ be a subgroup of the automorphism group of $X$.
	In \cite{Be-Me} J.Bertin, A. M\'ezard 
	proved that the equivariant cohomology of Grothendieck 
	$H^1(G,\T_X)$ is  the tangent space of the global deformation functor
	of smooth  curves with automorphisms. 
	The   dimension  of the $k$-vector space  $H^1(G,\T_X)$ is a measure of the  directions a curve can be deformed
	together with  a subgroup of the automorphism group.  The computation of this dimension turns out to be a difficult 
       problem, and  in the literature there are results concerning ordinary curves \cite{CK}, cyclic groups
       of order $p$ \cite{Be-Me}. In \cite{Ko:2002a} the author   attempts to compute   this space using the
        low terms   of  the Lyndon-Hochschild-Serre spectral sequence.
This allows us to compute the dimension of $H^1(G,\T_X)$ if $G$ is an elementary abelian group {\em i.e.} isomorphic 
to a direct product of cyclic groups of order $p$.
  The disadvantage of this method is that it 
        involves  the computation of the transfer map that depends on the group structure of the  decomposition series of the decomposition 
groups $G(P)$ at the wild ramification points.

The author in \cite{KoJPAA06} proposed an alternative  method of approaching the problem of computing the dimension of $H^1(G,\T_X)$. 
Precisely   one can prove that for a $p$-group $G$ and for the sheaf $\Omega_X$ of differentials on $X$ the following equation holds:
\[
 H^1(G,\T_X)=H^0(X,\Omega_X^{\otimes 2})_G.
\]
This observation allows us to compute the desired dimension if we know the Galois module structure of the space 
of $2$-holomorphic differentials. The knowledge of the Galois module structure  of the space of holomorphic differentials  in
 positive characteristic is still  an open problem and there are only partial results in the case of tame covers
\cite{Kani:88},\cite{Nak:87i},\cite{Nak:85j}
	 and when   $G$ is a cyclic group of order equal to 
the characteristic \cite{Nak:86}.  The computation of the Galois module structure of holomorphic differentials 
in the tame case does not offer any new information for the computation of $ H^1(G,\T_X)$ since 
the classical methods from the theory of Riemann surfaces can be used in this case.

In \cite{KoJPAA06} the author was able to recover the result of J.Bertin, A. M\'ezard  \cite{Be-Me} concerning the dimension 
of $H^1(G,\T_X)$ in the case  $G=\mathbb{Z} /p\mathbb{Z}$ by applying the computation  of Nakajima \cite{Nak:86} concerning the 
$\mathbb{Z}/p\mathbb{Z}$-module structure of $\Omega_X^{\otimes 2}$.

In this article we use the results  of B. K\"ock \cite{Koeck:04}  on weakly ramified covers, {\em i.e.} on covers so that 
$G_2(P)=\{1\}$ for all wild ramified points, in order to determine $\dim_k H^1(G,\T_X)$. 
Since all ordinary curves are weakly ramified \cite[th. 2i]{Nak}  we recover the result of Cornelissen and Kato \cite{CK} for 
deformations of ordinary curves.

Many authors \cite{Subrao},\cite{Nak:85},\cite{Kani:88},\cite{borne04},\cite{Stalder04}, in order to study the $k[G]$-module structure of spaces $\Omega_X(-D)$ considered the $k[G]$-module 
structures of the spaces of semisimple and nilpotent differentials with respect to the Cartier operator.  In section \ref{pranksection}
we follow this approach. The representation module coresponding to semisimple differentials is studied and the problem 
of computing the space of covariant differentials is reduced to the problem of computing nilpotent covariant differentials. More precisely we are able to 
prove that:
\begin{equation} 
 \dim_k H^1(G,\mathcal{T}_X)= 2(g_Y-1)+ \gamma_Y-1 +r + \dim_k  \Omega_X(-D)^n _G,
\end{equation}
where $g_Y,\gamma_Y$ denote the genus and  the $p$-rank of the Jacobian of $Y=X/G$ and $r$ is the number of 
points of $Y$ ramified in $X\rightarrow X/G$.

In his PhD thesis \cite{bornephd} N. Borne  proved an equivariant Riemann-Roch theorem and provided an equivariant 
Euler characteristic $K_0(G,X) \rightarrow R_k(G)$, where $K_0(G,X)$ is the Grothendieck group of $G-\mathcal{O}_X$-modules
and $R_k(G)$ is the Brauer Character Group. In the case of $p$-groups  in positive characteristic this approach does 
not give  us more  information than the classical Riemann-Roch theorem; the group $R_k(G)=\mathbb{Z}$ \cite[15.6]{SerreLinear} and equivariant information 
is lost.

N. Borne in \cite{Borne06}
proposed a refinement of his equivariant $K$-theory by using the ideas of Auslander in 
modular representation theory.  This new theory enables him to compute the Galois module structure of 
spaces of global sections of linear systems of curves of the form 
 $\mathcal{O}_X(D)$ where $\deg D \geq 2g_X -2$ and $G$ is a 
cyclic group of the form $\mathbb{Z}/p^v \mathbb{Z}$. We apply the results of Borne in order to 
prove:
\begin{pro} \label{main-pro}
Assume that the group $G=\mathbb{Z}/p^v \mathbb{Z}$ acts on $X$. Let $X_{ram}$ denote the set of ramification 
points of the cover $\pi:X \rightarrow X/G=Y$. For every point $P \in X_{ram}$ let $N_P$ denote the highest jump in the lower 
ramification filtration and $f(P)$ the highest jump in the upper ramification filtration. Let $\Delta:\mathbb{Z} \rightarrow \mathbb{Z}$
denote the map $\Delta:a\mapsto \lf a/p \rf$, and $\Delta^k$ the composition of $k$ times the map $\Delta$. Let $g_Y$ be the genus of the 
quotient curve $Y$ and let $r= \#\pi_*( X_{ram})$.

 The dimension of $H^1(\mathbb{Z}/p^v \mathbb{Z},\T_X)$ is given by:
\[
\dim_k H^1(\mathbb{Z}/p^v \mathbb{Z},\T_X)=3(g_Y-1)+  2r+\sum_{P \in \pi_*( X_{ram})} \left( 2 f(P) + \Delta^k (-2 N_P-2)\right).
\]
\end{pro}
%

This result can be obtained also by using Grothendieck's equivariant theory. 
One can   compute the  the group $H^0\big(Y, \pi_*^G (\T_X)\big)$  \cite{Be-Me},\cite[sec 3.]{Ko:2002a} 
and the local contributions at each wild ramified point given in \cite[prop.  4.1.1]{Be-Me}.

In deformation theory it is an interesting problem to compare the dimension of the tangent space of  the deformation functor 
to the dimension of the versal deformation ring. In the case of ordinary curves and in   the $p$-cyclic case  the Krull dimension 
has been computed \cite{CK},\cite{Be-Me}. 
In \cite{Pries:05} R. Pries studied unobstructed deformations of curves acted on by  Abelian groups  
under the additional assumption that the deformations do not split the branch locus.
 A comparison of the result of Pries to the result of proposition \ref{main-pro} gives
the dimension of the space generated by obstructed deformations and of deformations that split the branch locus.

{\bf Acknowledgements} The author would like  to thank G. Cornelissen  and  M. Matignon for
pointing him to the work of N. Borne, B. K\"ock, N. Stalder. 
The author would like to thank B. K\"ock for his remarks and corrections on an earlier version of 
this article.


\section{Spaces of 2-holomorphic differentials.}

We will denote by $\Omega_X$ the sheaf of differentials of $X$.
Let $D$ be an effective divisor on the curve $X$ and let $G$ be a finite subgroup of $\A (X)$.
The set of meromorphic  differentials on $X$ is denoted by $\mathcal{M}_X=\{xdy: x,y \in K(X)\}$, 
where $K(X)$ denotes the function field of the curve  $X$.
We will denote by 
\[
\Omega_X(D)=\{ \omega \in \mathcal{M}_X:\;\mathrm{div}(\omega) \geq D\}=H^0(X,\Omega_X \otimes \mathcal{O}_X(-D)).
\]
and by
\[
 \mathcal{L}_X(D)=\{f \in K(X): \mathrm{div}(f) \geq -D \}=H^0(X,\mathcal{O}_X(D)).
\]
We will denote the $k$-dimension of $\mathcal{L}_X(D)$ by $\ell(D)$.
The reader should be carefull about the notation here. We follow the notation of   Serre, used also in the papers of  Stichtenoth, Nakajima, Subrao. 
In the papers of N. Borne \cite{Borne06}, E. K\"ock \cite{Koeck:04}, E. Kani \cite{Kani:88}, N. Stalder 
\cite{Stalder04}  the notation  $\Omega_X(D)$ is used for what we write as 
$\Omega_X(-D)$. 
\begin{lemma} \label{efecDG}
	Let $G$ be a $p$-group and consider the cover  $X \rightarrow X/G$. We  assume that
the genus of the curve $X$ is $g_X\geq 2$.
 There is a $G$-invariant  differential $\omega$ in $X$, 
	such that $\mathrm{div}(\omega)$ is effective.
\end{lemma}
\begin{proof}
	Let $b_1,\ldots,b_r$ be the ramification points of the cover 
	$\pi:X \rightarrow Y=X/G$. 
In order to find $\omega$ we have to select a meromorphic differential $\phi$ of $Y$ and take the 
pullback $\pi^*(\phi)$. Then $\pi^*(\phi)$ is $G$-invariant and we have to select it so that 
$\mathrm{div}(\phi) \geq 0$.  Fix a meromorphic differential $\phi_1$ on $Y$ and 
consider the set $f \phi_1$ where $f$ is an arbitrary element in the function field of $Y$.
We have that 
\begin{equation} \label{upper-cond}
 \mathrm{div}( \pi^* (f \phi_1))=\pi^* \mathrm{div}(f\phi_1) +R \geq 0,
\end{equation}
	where $R$ is the ramification divisor given by 
	\[
	R=\sum_{i=1}^r  \sum_{P \mapsto b_i} 
	\sum_{\nu=0}^{\infty} \big({e_\nu(P) -1}\big) P.
	\]
In order to check that the divisor on the left hand side of (\ref{upper-cond}) is positive we push
forward  again and arrive at
\begin{equation} \label{down-cond}
 \mathrm{div}( f\phi_1) + \sum_{i=1}^r \sum_{\nu=0}^{\infty}  \frac{ e_\nu(b_i)-1}{e_0(b_i)} \geq 0.
\end{equation}
The later condition is equivalent to $f \in \mathcal{L}_Y(K+A)$,
	where $A$ is the effective divisor
	\[
	A:= \sum_{i=1}^r \lf \sum_{\nu=0}^{\infty} \frac{e_\nu(b_i) -1}{e_0(b_i)}
	\rf b_i.
	\]
Using Riemann-Roch on the curve $Y$ we compute 
	\[
	\ell(K+A)= \ell(-A) + 2g_Y-2 + \deg(A)-g_Y +1=g_Y -1 +\deg(A) \geq 1.
	\]
Indeed, 
since $G$ is a $p$-group and $g_X>2$ there should be at least one wild ramified point. Every ramified point contributes at least 
$1$ to the degree of $A$, since $e_0(b_i)=e_1(b_i)$ and if $g_Y \geq 1$ we are done. In the case   $g_Y=0$ we use  $g_X \geq 2$ 
and the 
Riemann-Hurwitz formula to see  that there should be either at least  two ramified points so    $\deg A \geq 2$, or one ramified point 
with $\sum_{\nu=0}^{\infty} \frac{e_\nu(b_i) -1}{e_0(b_i)}>2$, and in this case also $\deg A \geq 2$.
	
Every   $f\in \mathcal{L}_Y(K+A)$ gives rise to the desired differential $\omega=\pi^*(f \phi_1)$.  
	Moreover we 
	can select $f \in \mathcal{L}_Y(K +A)$ such that is has polar divisor $A$. 
	This imply that the support of $\pi^* (f \phi_1)$ has no intersection 
	with the branch locus. 
\end{proof}

\begin{pro}
Assume that the $p$-group  $G$ acts on the group $G$. There is an effective and invariant divisor 
$D^*=\mathrm{div}( p^*(\phi) )=\pi^* \big(\mathrm{div} (\phi )\big) +R \geq 0$ so that 
 the module  $\Omega_X(-D^*)$ is 
isomorphic to $H^0(X,\Omega_X^{\otimes 2})$ as a $k[G]$-module.
\end{pro}
\begin{proof}
According to lemma \ref{efecDG} we can select $\phi$ so that $\pi^*(\phi)$ is $G$-invariant and the divisor 
$D^*:=\mathrm{div}( \pi^*(\phi))=\pi^* \big(\mathrm{div} (\phi )\big) +R$ is an effective $G$-invariant divisor.
 Every differential can be written as $f \phi$ and every $2$-differential is 
an expression of the form $\omega \otimes \phi$, for an other differential $\omega$ on $X$. 
The space 
\begin{eqnarray*}
 H^0(X,\Omega_X^{\otimes 2})= & \{ \omega \otimes \phi : \mathrm{div} ( \omega \otimes \phi) \geq 0\}= \\
                                               = &  \{ \omega \in \Omega_X: \mathrm{div} ( \omega) \geq -\mathrm{div}(\phi) \}.=\\
                                                =&  \Omega_X(-D^*).
\end{eqnarray*}
\end{proof}

\section{The weakly ramified case}
In this section we assume that the cover $X \rightarrow X/G$ is {\em weekly} ramified, {\em i.e.} for 
every $P \in S$ we have that $G_i(P)=\{1\}$ for all $i \geq 2$. The group $G$ is always assumed to be a $p$-group.
The ramification divisor is computed  $R=\sum_{P \in S}  2(e_0(P)-1)P$.

\begin{lemma}
 The module $\Omega_X(-D^*  -\sum_{P\in S} 3P)$ is a projective $k[G]$-module.
\end{lemma}
\begin{proof}
Let $K$ be a canonical divisor. 
Let us write 
\[\Omega_X(-D^*  -\sum_{P\in S} 3P)=\mathcal{L}_X( K + D +\sum_{P \in S} 3P).\]
Let $D':=K+D^* + \sum_{P \in S} 3P=\sum n_P P$. 
We compute that 
\[
 D'=K+D^* + \sum_{P \in S} 3P = 2\pi^* (\mathrm{div}(\phi)) + \sum_{P\in S} 4(e_0(P)-1) +3)P.
\]
Since we have assumed that $G$ is a $p$-group we have that $e_0(P)=e_1(P)=e^w(P)$. Therefore 
$n_P \equiv -1 \mod e_0$ for all wild ramified points of $X \rightarrow X/G$.
On the other hand 
\[
 H^1(X,\O_X(D'))=\mathcal{L}_X(K-D')=0,
\]
since $\deg(K-D')<0$.
The desired result follows by  \cite[th. 2.1]{Koeck:04}.

B. K\"ock proposed to me the following more abstract  approach: By \cite[th. 44]{bornephd} there is a $G$-invariant canonical divisor $K_X=\sum_{P} m_P P$ on $X$. 
The proof of corollary 2.3 in \cite{Koeck:04} implies that all $m_P \equiv -2 \mod e_0^w(P)$. Therefore the divisor,  $2K_X + \sum_{P \in S} 3P
=\sum n_P P$
has $n_P \equiv 2 (-2)+3=-1 \mod e_0^w(P)$ as required.
\end{proof}

We can now form the short exact sequence:
\begin{equation} \label{ses1}
 0 \rightarrow \Omega_X(-D^*) \rightarrow \Omega_X(-D^* - \sum_{P \in S} 3P) \rightarrow \frac{\Omega_X(-D^*  - \sum_{P \in S} 3P)}{\Omega_X(-D^*)} 
\rightarrow 0.
\end{equation}
The short exact sequence of sheaves 
\[
 0 \rightarrow \Omega_X\otimes \mathcal{O}_X(-D^*) \rightarrow \Omega_X\otimes \mathcal{O}_X (-D^* - \sum_{P \in S} 3P) \rightarrow 
\Sigma\rightarrow 0
\]
where
\[
 \Sigma:=\frac{\Omega_X\otimes \mathcal{O}_X (-D^* - \sum_{P \in S} 3P)}{\Omega_X\otimes \mathcal{O}_X(-D^*)}
\]
gives rise to a long exact sequence of $k$-vector spaces by applying the functor of global sections. This long exact sequence combined with 
  $H^1(X,\Omega_X\otimes \mathcal{O}_X(-D^*))=0$ allows us to express  
\[
H^0(X,\Sigma) =\frac{\Omega_X(-D^*  - \sum_{P \in S} 3P)}{\Omega_X(-D^*)}\]
as the direct sum of the stalks $\Sigma_P$ of the skyscrapper sheaf $\Sigma$ at points $P\in S$.
Let us denote by $\Sigma':=H^0(X,\Sigma)$.

Equation (\ref{ses1}) gives the following long exact sequence:
\[
 0  \rightarrow H_1(G, \Sigma') \rightarrow \Omega_{X} (-D^*)_G \rightarrow
\Omega_{X}(-D^* - \sum_{P \in S} 3P)  \rightarrow \Sigma'_G \rightarrow 0, 
\]
since $\Omega_X(-D^* - \sum_{P \in S} 3P)$ is a $k[G]$-projective module.
Therefore, the desired dimension can be computed:
\begin{equation} \label{compOrd}
\dim_k\Omega_X(-D^*)_G = \dim_kH_1(G,\Sigma') + \dim_k\Omega_{X} \big(-D^*  - \sum_{P \in S} 3P \big)_G - \dim_k\Sigma'_G 
\end{equation}
In what follows we will compute every summand on the right hand side of (\ref{compOrd}).

\begin{lemma} \label{projcomp}
 $\dim_k \Omega_{X} \big(-D^*  - \sum_{P \in S} 3P \big)_G=3 (g_Y -1) + 3 r.$
\end{lemma}
\begin{proof}
 Using the theorem of Riemann-Roch we compute:
%
%
\begin{equation} \label{RR1}
 \dim_k \Omega_{X} \big(-D^*  - \sum_{P \in S} 3P \big) =2g_X-2  + |G| (2g_Y-2)
                                                                                  +  \sum_{P\in S}  \big( 2(e_0(P)-1)+ 3\big) 
\end{equation}
\[
                                                                                  - g_X +1 + \dim_k \mathcal{L}_X(-D^* -\sum_{P \in S} 3P).
\]
But $\deg( -D^*  -\sum_{P\in S} 3P) < 0$ therefore $\dim_k \mathcal{L}_X(-D^* -\sum_{P \in S} 3P)=0$. 
Riemann-Hurwitz implies that 
\begin{equation} \label{RH1}
 2g_X-2 = | G| (2 g_Y -2) + 2 \sum_{P\in S} (e_0(P)-1).
\end{equation}
By combining (\ref{RR1}),(\ref{RH1}) we obtain:
\[
 \dim_k \Omega_{X} \big(-D^*  - \sum_{P \in S} 3P \big)=3 |G|(g_Y-1)  + 3 r|G|.
\]
Since $\Omega_{X} \big(-D^*  - \sum_{P \in S} 3P \big)$ is projective it is of the form $k[G]^a$,  where $a=3(g_Y-1) +3 r$,
and each $k[G]$ direct sumand contributes $1$ to the dimension of $\Omega_{X} \big(-D^*  - \sum_{P \in S} 3P \big)_G$. 
The desired result follows.
\end{proof}

Let us now study the space $\Sigma'$ as a $G$-module:
The differential $\pi^*(\phi)$ defined in lemma \ref{efecDG} is an invariant differential and we have that 
$\mathrm{div} (\pi^*(\phi))=D^*$. Every differential $\omega$ can be written as $\omega=f \pi^*(\phi)$,
Therefore, 
\[
 \Omega_X(-D^*)=\{ \mathrm{div} f \pi^*(\phi)  > -D^* \}\cong \mathcal{O}_X (2D^*),
\]
and 
\[
 \Omega_X(-D^* -\sum_{P\in S}3P)  \cong \mathcal{O}_X(2D^*+\sum_{P\in S} 3P),
\]
where the last two isomorphisms are isomorphisms of $k[G]$-modules.

Let $O(P)=\{ g P: g\in G\}$ denote the orbit of $P$ under the action of the group $G$.
The set $S$ of ramification points can be written as a disjoint union of orbits of points of $X$:
\[
 S =\bigcup_{j=1}^r O(P_j),
\]
for a selection $P_1,\ldots,P_r$ of points of $X$.
We can write $H^0(X,\Sigma)$ as 
\begin{equation} \label{sigma2}
\Sigma'=\bigoplus_{j=1}^r \bigoplus_{P \in O(P_j)} \Sigma_P,
\end{equation}
$\Sigma_P$ is the stalk of $\Sigma$ at $P$. 
Thus, (\ref{sigma2}) can be written as 
\[
 \Sigma'= \bigoplus_{j=1}^r  \mathrm{Ind}_{G(P_j)}^G \Sigma_{P_j}.
\]
Shapiro lemma \cite[6.3.2, p.171]{Weibel} implies that 
\begin{equation} \label{shapirosum}
 H_*(G,\Sigma')= \bigoplus_{j=1}^r H_*(G, \mathrm{Ind}_{G(P_j)}^G \Sigma_{P_j}) =\bigoplus_{j=1}^r H_*( G(P_j), \Sigma_{P_j}).
\end{equation}

Let $P$ be a ramified point.
We  have assumed that the cover is weakly ramified we have that $G_2(P)=\{1\}$. This implies that for a local 
uniformizer $t_P$ at $P$ we have 
\begin{equation} \label{action}
 G(P) \ni g: t_p \mapsto t_P \left (1 + \alpha_1(g)  t_P +\sum_{\nu=2}^\infty \alpha_\nu(g) t_P^\nu \right),
\end{equation}
where $\alpha_1(g) \neq 0$. The map $g\mapsto a_1(g)$ is a homomorphism and allows us to see  $G(P)$
as a finite dimensional $\mathbb{F}_p$-vector subspace of $k$.

The quotient $\Sigma_P$ is then generated as a $k$-vector space by elements 
\[
 \Sigma_P= \left\{ \omega_1:= \frac{1}{t_P^{4e(P)-1}}, \omega_2:=\frac{1}{t_P^{4e(P)-2}},\omega_3:= \frac{1}{t_P^{4e(P)-3}}\right\}
\]
and the action is given by:
\begin{equation} \label{action1}
 \sigma: \frac{1}{t_P^{4e(P)-1}} \mapsto \frac{1 +  \alpha_1 (g) t_P  + \alpha_2(g)  t_P^2}{t_P^{4e(P)-1}} =
\frac{1}{t_P^{4e(P)-1}}+ \alpha_1(g)  \frac{1}{t_P^{4e(P)-2}} + \alpha_2(g) \frac{1}{t_P^{4e(P)-3}},
\end{equation}
since $\binom{-4e(P)+1}{1}=1$ and $\binom{-4e(P)+1}{2}=0$.
\begin{equation} \label{action2}
 \sigma:\frac{1}{t_P^{4e(P)-2}} \mapsto \frac{1 +  2 \alpha_1 (g) t_P}{t_P^{4e(P)-2}}= \frac{1}{t_P^{4e(P)-2}}+2\alpha_1(g)  \frac{1}{t_P^{4e(P)-3}}
\end{equation}
and 
\begin{equation} \label{action3}
 \sigma:\frac{1}{t_P^{4e(P)-3}}\mapsto \frac{1}{t_P^{4e(P)-3}}.
\end{equation}

\begin{lemma} \label{compS1}
 We have that $\dim_k \Sigma'_G=H_0(G,\Sigma')=r$
\end{lemma}
\begin{proof}
Let $P$ be one of the $\{P_1,\ldots,P_r\}$.
Let $\omega \in \Sigma_{P}$. Using equations (\ref{action1}),(\ref{action2}),(\ref{action3}) 
we observe that there is only one invariant  element in $\Sigma_P$, 
therefore the images  of the linear maps  $(g-1)\omega$ are all isomorphic and two dimensional generated by the 
elements $\omega_2,\omega_3$ and so  we arrive at: 
\[
 \dim_k H_0(G(P),\Sigma_P)=\dim_k \frac{\Sigma_P}{g \omega -\omega}=1.
\]
For the global computation we use  (\ref{shapirosum}) in order to obtain:
\[
 \dim_k \Sigma'_G=H_0(G,\Sigma')=r.
\]
\end{proof}
\begin{lemma} \label{compS2}
 We have that $\dim_k H_1(G,\Sigma')=\sum_{j=1}^r \log_p |G(P_j)|-r$
\end{lemma}
\begin{proof}
 Let $P$ be one of the $\{P_1,\ldots,P_r\}$.
We will use the {\em normalised bar resolution} defined in \cite[6.5, p. 177]{Weibel}
in order to compute $H_1(G, \Sigma_P)$.
Recall that $k[G(P)] \cong B_0$ and $B_0$ is generated by the symbol $[\;]$, 
$B_1$ is the free $k[G(P)]$-module on the set of symbols $\{[g]: g\in  G(P) \backslash \{1\}\}$, 
and $B_2$ is the free $k[G(P)]$-module on the set of symbols $\{ [g|h]:g,h \in  G(P) \backslash \{1\} \}$.
For the differential maps  we have  $\partial_1:B_1 \rightarrow B_0$ and  $\partial_1 [g]=(g-1)[\;]$, and 
$\partial_2:B_2 \rightarrow B_1$, $\partial_2[g|h]=g[h]-[gh]+[g]$. Higher  $B_n$ can be similarly defined but we don't need them here.
The  group $H_1(G,\Sigma_P)$ is defined by the homology at position $1$ of the chain complex $\Sigma_P \otimes B_*$.  
We have:
\begin{equation} \label{chain1}
 H_1(G, \Sigma_P)=\frac{\ker(\partial_1:\Sigma_P \otimes B_1 \rightarrow \Sigma_P \otimes B_0)}{ \mathrm{im}(
\partial_2: \Sigma_P \otimes B_2 \rightarrow \Sigma_P \otimes B_1)}.
\end{equation}
We will focus first on the study of the local components $\Sigma_P$.
Let $\omega_3=\frac{1}{t_P^{4e(P)-3}}$, $\omega_2=\frac{1}{t_P^{4e(P)-2}}$, $e_1=\frac{1}{t_P^{4e(P)-1}}$.
Observe that  the  space generated by $\omega_3 \otimes [g]$ is in $\mathrm{ker} \partial _1$, since 
\[
 \partial _1 \big( \omega_3 \otimes [g] \big) = \omega_3 (g-1) \otimes [\;] =0.
\]
On the other hand 
\[
 \partial_2 \big( \omega_2 \otimes [g|h] \big)=\omega_2 \otimes( g[h]-[gh]+[g]) =\omega_2 g \otimes [h] -\omega_2 \otimes [gh] + 
\omega_2 \otimes [g]=
\]
\[
 =\alpha_1(g) \omega_3 \otimes [h] + \omega_2 \otimes ( g[h]-[gh]+[g]).
\]
Since $\omega_2,\omega_3$ are linear independent and $\alpha_1(g) \neq 0$, 
 we obtain that all elements of the form $\omega_3 \otimes [h]$, $h \in G(P)$ 
are $0$ in the homology group $H_1(G(P),\Sigma_P)$.

Observe that an element $\sum_{g\in G(P)} \lambda_g \omega_2 \otimes [g]$ is in the kernel of $\partial_1$ if and only 
if $\sum_{g \in G(P)}  \lambda_g \alpha_1(g)=0$, since we have assumed that  $p\neq 2$.
On the other hand 
\begin{equation} \label{impar}
 \partial_2 \big(  \omega_1 \otimes [g|h] \big)= \omega_1 \otimes (g[h] -[gh]+[g]) =\omega_2 \alpha_1(g) \otimes [h] +
\omega_1 ( [h]-[gh]+[g]),
\end{equation}
therefore $\omega_2 \otimes [h]$ is zero in the homology group $H_1(G,\Sigma_P)$.

Consider an element 
\[
\omega =\sum_{g\in G(P)} \lambda_g \omega_1 \otimes [g] +  \sum_{g\in G(P)} \mu_g \omega_2 \otimes [g]+\sum_{g\in G(P)} \nu_g \omega_3 \otimes [g],
\]
and the image $\partial_1 \omega$ given by 
\begin{equation} \label{parom}
 \partial_1 (\omega)=\sum_{g\in G(P)} \alpha_1(g) \lambda_g  \omega_2 \otimes [\;] + 
\sum_{g\in G(P)} \left( 2 \alpha_2(g) \lambda_g + \mu_g \alpha_1(g) \right) \omega_3 \otimes [\;]. 
\end{equation}
Equation (\ref{parom}) gives two necessary conditions for $\omega \in \mathrm{ker} \partial_1$, namely:
\[
\sum_{g\in G(P)} \left( 2 \alpha_2(g) \lambda_g + \mu_g \alpha_1(g) \right)=0
\]
and 
\[
 \sum_{g\in G(P)} \alpha_1(g) \lambda_g=0.
\]
We have seen that $H_1(G(P),\Sigma_P)$ is generated by the images of elements of the 
form $\omega_1 \otimes [g]$. On the other hand, equation (\ref{impar}) implies that:  
\begin{equation} \label{chain34}
\omega_1 \otimes [gh] =\omega_1 \otimes [g] + \omega_1 \otimes [h] \mbox{ in } H_1(G(P),\Sigma_P).
\end{equation}
The groups $G(P)$ are elementary Abelian and can be written as  $G(P)=\bigoplus_{i=1}^{t} g_i \mathbb{F}_p$.
Equation (\ref{chain34}) implies that 
\[
 H_1(G, \Sigma_P) \subset   \omega_1  \otimes  \langle  [g_i] : i=1,\ldots, t \rangle_k.
\]
For a linear combination $\sum_{i=1}^t \lambda_i \omega_1  \otimes [g_i]$ we have that 
it is in $\ker \partial_1$ if and only if $\sum_{i=1}^t \lambda_i=0$. Therefore:
\[
 \dim_k H_1(G,\Sigma_P)=\log_p |G(P)|-1.
\]
The global contribution is computed using (\ref{shapirosum}) and equals to 
\[
 H_1(G, \Sigma')=\sum_{j=1}^r \log_p |G(P_j)|-r,
\]
{\em i.e.,} the desired result.
\end{proof}

Observe that the dimension of the tangent space of the deformation functor  can now be computed by combining  (\ref{compOrd}) and lemmata 
\ref{projcomp},\ref{compS1},\ref{compS2}:
\[
 \dim_k H^1(G,\T_X)= 3g_Y -3 +r + \sum_{j=1}^r \log_p |G(P_j)|. 
\]
Notice, that if the curve $X$ is ordinary and $G=\A(X)$ then the cover $X \rightarrow X /G$ is weakly ramified \cite[th. 2i]{Nak}  and
the above result coincides with the result of  G. Cornelissen and F. Kato \cite{CK} on deformations of ordinary curves.

\section{The $p$-rank representation.} \label{pranksection}

Let $D$ be an effective divisor on a curve $X$. 
On the spaces $ \Omega_X(-D))$
 one can define the action of the Cartier operator. For an introduction to all necessary 
material the interested reader may consult \cite{SerreMexico},\cite{Subrao},\cite{Nak:85},\cite{Stalder04}. There is the following decomposition of the above space
\[
  \Omega_X(-D) = \Omega_X(-D)^s \bigoplus  \Omega_X(-D)^n 
\]
where $\Omega_X(-D)^s, \Omega_X(-D)^n$ are the spaces of semisimple 
and nilpotent differentials  with respect to the Cartier operator.
The above decomposition is compatible with the $G$-action. For the  $k[G]$-module $\Omega_X(-D)^n$ of 
nilpotent elements little seems to be known. On the other hand the $k[G]$-module 
$V_D:=\Omega_X(-D)^s$ was studied by many authors (\cite{Nak:85},\cite{Kani:88},\cite{borne04},\cite{Stalder04}).
The $k[G]$-module $V_D$ is called  in the literature  {\em the $p$-rank representation}. 
In general, we have  the following  decomposition:
\[
 V_D=\mathrm{core}(V_D) \bigoplus  \bigoplus_{S \in \mathrm{Irr} G} P_G(S)^{b(G,D,S)},
\]
where $S$ runs over the set of equivalent classes of irreducible representations, $P_G(S)$ denotes the projective 
cover of $S$, and $b(G,D,S)\in \mathbb{N}$ are called {\em the Borne invariants}  corresponding to $G,D,S$ \cite{Stalder04}.

Since $G$ is a $p$-group  and $k$ is assumed to be of characteristic $p>0$ the only irreducible representation 
is the trivial one and has projective cover $k[G]$ \cite[15.6]{SerreLinear}. Moreover $\mathrm{core}(V_{D^*})=0$ since $D^*$ is non 
empty and contains all ramification points \cite{Nak:85}, \cite[4.5]{Stalder04}. Therefore, 
\[
 V_{D^*}=k[G]^{B(G,D^*,k)},
\]
 where $B(G,D^*,k)$ is an integer. 
\begin{pro}
 With the above notation $B(G,D^*,k)=2(g_Y-1)+ \gamma_Y-1 +r$, where $g_Y,\gamma_Y$ are the genus and the $p$-rank 
of the Jacobian of the curve $Y=X/G$, $r$ is the number of points that are ramified in the cover $X \rightarrow Y=X/G$.
\end{pro}
\begin{proof}
 Let us denote by $D^*_{red}$ the divisor that has the same support with $D^*$, so that for all prime divisors $P$ we 
have $v_P(D^*_{red})>0 \Rightarrow v_P(D^*_{red})=1$. According to  \cite[p. 175]{Subrao} 
\begin{equation} \label{subrao-Computation}
 \dim_k  \Omega_X(-D)^s=\dim_k \Omega_X(-D^*_{red})^s=\dim_k \Omega_X(0)^s+\deg D^*_{red}-1.
\end{equation}
But the space $\Omega_X(0)^s$ of semisimple regular differentials is of dimension $\gamma_X$, 
and the degree of $D^*_{red}$ is equal to $2(g_Y-1)|G| + r$. 

The Deuring-Shafarevich formula (\cite{Madan77},\cite{Subrao},\cite{Deuring},\cite{Shaf47},\cite{Nak:85}) relates the $p$-ranks 
$\gamma_X,\gamma_Y$:
\[
 \gamma_X-1=|G| \left( \gamma_Y-1 + \sum_{i=1}^r \left(1 - \frac{1}{{e_0(b_i)}} \right) \right).
\]
This combined with (\ref{subrao-Computation}) gives us that 
\[
 \dim_k V_{D^*}= \dim_k \Omega_X(-D^*)^s=|G| (\gamma_Y -1 + 2(g_Y-2) +r),
\]
and since $ V_{D^*}$ is projective we have that $V_{D^*}= k[G]^{\gamma_Y -1 + 2(g_Y-2) +r}$,
and the desired result follows. 
\end{proof}
The above computation allows us to compute the dimension of the space of nilpotent elements.
Indeed, the dimension of the space $\Omega_X(-D^*)$ is computed to be equal to 
\[
 \Omega_X(-D^*)=(3g_Y-3)|G| + |G| \sum_{i=1}^r \sum_{\nu=0}^\infty \frac{e_\nu(b_i)-1}{e_0(b_i)}. 
\]
Therefore, we have that 
\[
 \dim_k \Omega(-D^*)^n=|G|(g_Y -\gamma_Y) + |G| \sum_{i=1}^r \sum_{\nu=0}^\infty \frac{e_\nu(b_i)-1}{e_0} -r.
\]
If the curve $X$ is ordinary {\em i.e. $g_Y=\gamma_Y$} then $e_\nu(b_i)=0$ for all $\nu \geq 2$ and 
the above formula gives us 
\[
 \dim_k \Omega_X(-D^*)^n = |G| \left(\sum_{i=1}^r (1- \frac{2}{e_0(b_i)}) \right).
\]
Since this dimension is not divisible by $|G|$ the module $\Omega_X(-D^*)^n$ can not be projective. 

The conclusion concerning the dimension is that
\begin{equation} 
 \dim_k H^1(G,\mathcal{T}_X)= 2(g_Y-1)+ \gamma_Y-1 +r + \dim_k  \Omega_X(-D^*)^n _G.
\end{equation}

\section{Borne Theory} \label{sectionBorne}

Let $V_j$ denote the indecomposable $k[G]$-module of dimension $j$.
Denote by $V$ the $k[G]$-module  with $k$-basis $\{e_1,\ldots,e_{p^v}\}$ and action 
given by $\sigma e_\ell=e_\ell+e_{\ell-1}$, $e_0=0$.
Then, $V_j$ is the subspace of $V$ generated by $\{e_1,\ldots, e_j\}$.

Following Borne
we define:
\begin{definition}
Let $\pi:X \rightarrow Y$ be a Galois cover of curves defined over $k$ with Galois group $G\cong \mathbb{Z}/p\mathbb{Z}$.
For a ramified point $P$ of $X$ we define $N_P$ so that $\sigma (t_P)-t_P$ has valuation $N_P+1$, where $t_P$ is a local 
uniformizer at $P$. Let $X_{ram}$ denote the set of ramification locus of the above cover. 
For every $0 \leq \alpha \leq p-1$ we define a map $\pi_*^{\alpha}:\mathrm{Div}_X \rightarrow \mathrm{Div}_Y$ by 
\[
 \pi_*^\alpha D= \lf \frac{1}{p} \pi_*\left( D -\alpha \sum_{P \in X_{ram}} N_P P \right)  \rf,
\]
where $[\cdot]$ denotes the integral part of a divisor, taken coefficient by coefficient.
\end{definition}
We will use the following 
\begin{theorem} \label{BorneCyclic}
 Suppose that $X$ is acted on faithfully by the cyclic $p$-group $G\cong \mathbb{Z}/p^v \mathbb{Z}$.
We break the cyclic extension $X \rightarrow X/G$ to a sequence of cyclic $p$-extensions by defining for every 
$1 \leq n \leq v$ the cover $X_n =X/H_n$ where $H_n$ is the unique subgroup of $G$ of order $n$. We set $X_0=X$.
Let $\pi_n: X_{n-1} \rightarrow X_n$ denote the canonical morphism.
Let $D$ be a $G$-invariant divisor on $X$, then $H^0(X,\mathcal{O}_X(D)) \cong \bigoplus_{j=1}^v V_j^{\oplus m_j}$.
Suppose that $\deg(D) \geq 2g_X -2$. Then, the integers $m_j$ are given by 
\[
 m_j=\deg( \pi_{v*}^{\alpha_0(j)} \cdots \pi_{1*}^{\alpha_{v-1}(j)} D) - \deg( \pi_{v*}^{\alpha_0(j+1)} \cdots \pi_{1*}^{\alpha_{v-1}(j+1)} D) \mbox{ if } 1 \leq j \leq p^v-1,
\]
\[
 m_{p^v}= 1 -g_{X_v} + \deg(\pi_{v*}^{p-1} \cdots \pi_{1*}^{p-1} D),
\]
where for $1\leq j \leq p^v$ the integers $\alpha_0(j),\ldots,\alpha_{v-1}(j)$ are the digits of the $p$-adic expansion of $j-1$ {\em i.e.},
\[
 j-1=\sum_{h=0}^{v-1} \alpha_h(j)p^h, \mbox{ with } 0 \leq \alpha_h(j) \leq p-1.
\]
\end{theorem}
\begin{proof}
 \cite[th. 7.25]{Borne06}
\end{proof}

We would like to consider the space $\Omega_X(-D^*)=\mathcal{O}_X(K+D^*)=\mathcal{O}_X(2D^*)$, where 
$D^*=\pi^* (\mathrm{div}(\phi)) + \sum_{P\in X_{ram}}   \sum_{i=0}^\infty  (e_i(P) -1)P$. Notice that $\deg(2D^*)=4g_X-4 \geq 2g_X-2$ 
so theorem \ref{BorneCyclic} is applicable.

Define 
\[
 D(j):=\pi_{v*}^{\alpha_0(j)} \cdots \pi_{1*}^{\alpha_{v-1}(j)} 2D^*.
\]
We have  $m_j=\deg(D(j)) - \deg(D(j+1))$ for $1 \leq j \leq p^v-1$. 
Since  $\dim_k H_0(G,V_j)= \dim_k(V_j)_G=1$ the desired dimension is given by:
\begin{equation} \label{des12}
\dim_k \mathcal{O}_X(2D^*)_G=\sum_{j=1}^{p^v}m_j= \deg(D(1)) -\deg(D(p^v)) + 1-g_X + \deg(\pi_{v*}^{p-1} \cdots \pi_{1*}^{p-1} 2D^*).
\end{equation}
Observe that the $p$-adic expansion of $p^v-1$ is
\[
 p^v-1 =p-1 + (p-1)  p + \cdots (p-1) p^{v-1},
\]
so  $\alpha_i(0)=0, \alpha_i(p^v-1)=p-1$ for all  $0\leq i \leq v-1$.
Therefore (\ref{des12}) is 
simplified to 
\begin{equation} \label{des13}
 \dim_k \mathcal{O}_X(2D^*)_G=\deg(D(1))+ 1-g_X.
\end{equation}

\begin{definition}
 Let $X_{ram,i}$ be the set of points of $X_i$  ramified in the cover $X_i \rightarrow X_{i+1}$.
\end{definition}
Notice that $\mathbb{Z}/p\mathbb{Z}=Gal(X_0/X_1)$ is contained in every subgroup of $G$, 
therefore $X_{ram}=X_{ram,0}$.

Denote by $\Delta$ the map $\mathbb{Z} \rightarrow \mathbb{Z}$ sending $\Delta: a \mapsto \lf \frac{1}{p} a\rf$, 
and by $\Delta^n$ the composition $\Delta \circ \cdots \circ \Delta$
 of $n$ times $\Delta$. If $a=\sum_{i=0}^\ell \alpha_i p^i$ is the $p$-adic expansion of $a$, then 
\[
 \Delta^n(a)=\sum_{i=n}^\ell \alpha _i p^{i-n}.
\]
Assume that $H=\mathbb{Z}/p\mathbb{Z}$.
Let $O(P)=\{gP :g \in H\backslash H(P)\}$ be the orbit of $P$ under the action of $H$, and consider the $H$-invariant divisor:
$
 \sum_{Q \in O(P)}  a Q.
$
We have that 
\[
 \pi_*^{1} \left( \sum_{Q \in O(P)}  a Q \right) =
\left\{
\begin{array}{ll}
a \cdot \pi_*(Q) & \mbox{ if } Q \mbox{ does not ramify in } X \mapsto X/G\\
\Delta(a) \cdot  \pi_*(Q) & \mbox{ if } Q \mbox{ ramifies in }  X \mapsto X/G
 \end{array}
\right.
\]
The divisor $2D$ can be written as 
\[
 2 \pi^*\mathrm{div}(\phi) +  \sum_{P \in X_{ram}} \sum_{i=0}^\infty (2e_i (P) -2)P.
\]
Since the divisor $\pi^*\mathrm{div}(\phi)$ has empty intersection with the ramification locus, we have
\[
 \pi_{*}^{1}(2D^*)= (\pi_v \circ \pi_{v-1} \circ \cdots \circ \pi_2)^* 2 \mathrm{div}(\phi) +\sum_{P \in \pi_{1*}( X_{ram}) }  \lf \frac{1}{p}\sum_{i=0}^{N_P} (2e_i(P)-2)\rf P=
\]
\[
 (\pi_v \circ \pi_{v-1} \circ \cdots \circ \pi_2)^* 2 \mathrm{div}(\phi) +
 \sum_{P \in \pi_{1*}(X_{ram}) } \Delta \left( \sum_{i=0}^{N_P} (2e_i(P)-2) \right)P.
\]
The decomposition group $G_0(P)$ is a cyclic subgroup of the whole group $G$, therefore
$G_0(P)=\mathbb{Z} /p^{k(P)} \mathbb{Z}$. Observe that a point $P \in X_{ram}$ is  fully ramified in 
all covers $X_{i} \rightarrow X_{i+1}$ for $i \leq k(P)$.
We can see that 
\[
 D(1)=2\mathrm{div}(\phi) + \sum_{P \in \pi_* (X_{ram})} \sum_{i=1}^{k(P)} \Delta^{k(P)} \left(  \sum_{i=0}^{N_P} (2e_i(P)-2) \right) P. 
\]
Let 
\[
 1 \leq i_1 \leq i_2 \leq \cdots \leq i_{k(P)}=N_P, 
\]
be the sequence of the jumps in the ramification filtration at the point $P$, {\em i.e.}
\[
 G_{i_k} \varsupsetneq G_{i_k+1} \mbox{ and } G_{i_{k(P)}+1}=\{1\}.
\]
\begin{lemma} \label{cyHasArf}
There are strictly positive integers $a_0,a_1,\ldots,a_{k-1}$, so that 
the sequence of jumps for the ramification filtration for the cyclic group $\mathbb{Z}/p^k \mathbb{Z}$ is of the form:
\[
  i_\nu=\sum_{\mu=0}^{\nu-1} a_\mu p^\mu.
\]
In particular, for all $\nu \geq 2$
\[
 i_\nu-i_{\nu-1} =a_{\nu-1}  p^{\nu-1} 
\]
\end{lemma}
\begin{proof}
 This is a direct consequence of the Hasse-Arf theorem for Abelian groups and it is explained in the example that appears on  page 76 in \cite{SeL}.
\end{proof}

Notice that we have $k(P)$ jumps in the ramification filtration since $G_{i_k}/G_{i_{k}+1}$ is 
elementary Abelian, therefore isomorphic to $\mathbb{Z}/p\mathbb{Z}$.
We set $d_P=\sum_{i=0}^{N_P} (2e_i(P)-2)$, and we compute:
\[
 d_p= \sum_{i=0}^{i_1} (2e_{i_1}(P)-2) + \sum_{i=i_1+1}^{i_2}(2e_{i_2}(P)-2) + \cdots + 
\sum_{i_{k(P)-1}+1}^{i_{k(P)}} (2e_{i_k(P)}(P)-2).
\]
Since $e_{i_\ell}=p^{k(P)-\ell+1}$ for all $1 \leq \ell \leq k(P)$ we have that 
\[
d_P=2 (p^{k(P)}-1) (i_1+1) + 2((p^{k(P)-1}-1)(i_2 -i_1) + \cdots + 2(p-1)( i_{k(P)} - i_{k(P)-1})=
\]
\[
 =\sum_{\nu=1} ^{k(P)} 2 (p^{k(P)-\nu+1}-1)(i_\nu -i_{\nu-1}), \mbox{ where } i_0=-1.
\]
Using lemma \ref{cyHasArf} we obtain that there are integers $a_0(P),a_1(P),\ldots,a_{k(P)-1}(P)$ so that 
\begin{equation} \label{cycd}
 d_p= 2 p^{k(P)}-2+\sum_{\nu=1}^{k(P)} 2p^{k(P)-\nu+1} p^{\nu-1} a_{\nu-1} (P) +     \sum_{\nu=1}^{k(P)}  2(-i_\nu +i_{\nu-1}) =
\end{equation}
\[
2 p^{k(P)}+p^{k(P)} 2 \sum_{\nu=1}^{k(P)}  a_{\nu-1}   -2 N_P-2
\]
One can compute \cite[exam. p.72 ]{SeL} that the sum $f(P):=\sum_{\nu=1}^{k(P)}a_{\nu-1}(P)$ is the highest jump in the upper 
ramification filtration. 
Using (\ref{cycd}) we compute that 
\[
 \Delta^{k(P)}(d_p)= 2 f(P) + \Delta^k (-2 N_P-2).
\]

Combining all the above we arrive at 
\begin{equation} \label{findegd1}
 \deg D(1) =4g_Y-4 + 2\#\pi_*( X_{ram})+  \sum_{P \in \pi_*( X_{ram})}  \left( 2 f(P) + \Delta^k (-2 N_P-2)\right)
\end{equation}

The desired result follows:
\[
 \dim_k\mathcal{O}(2D^*)=3(g_Y-1)+  2\#\pi_*( X_{ram})+\sum_{P \in \pi_*( X_{ram})} \left( 2 f(P) + \Delta^k (-2 N_P-2)\right).
\]
Observe that if $G =\mathbb{Z}/p\mathbb{Z}$ the above result coinsides with the computation of \cite{Be-Me},\cite{KoJPAA06}.

{\bf Example:}
Let $\sigma_i=\sum_{\nu=0}^{i-1} a_\nu$ $i=0,\ldots,k(P)$ be the jumps in the upper ramification filtration at the 
wild ramified point $P$.
It is known that $\sigma_{i+1}=p \sigma_i$ or $\sigma_{i+1}/p \geq \sigma_i$ and $p\nmid \sigma_{i+1}$ \cite{PriesNumTh},\cite{Schmid36}.
Assume that all jumps in the ramification filtration are of the form $\sigma_{i+1}=p \sigma_i$. In this case one can prove by induction 
that $a_i=a_0 p^{i-1} (p-1)$, and we can compute 
\[
 f(P)+ \Delta^k (-2 N_P-2)=\sum_{\nu=0}^{k(P)-1} a_\nu + \lf \frac{-2 \sum_{\nu=0}^{k(P)-1} a_\nu  p^\nu-2}{p^{k(P)}} \rf= 
\]
\[
 = \sum_{\nu=0}^{k(P)-1} a_0 (p-1) p^{\nu-1} +  \lf \frac{-2 \sum_{\nu=0}^{k(P)-1} (p-1)  p^{2\nu-1}-2}{p^{k(P)}} \rf
\]
\[
 =a_0 (p^{k(P)}-1) + \lf
-a_0 2\frac{
p^{2k(P)}-1
}{
p^{k(P)+1}(p+1)
}          -\frac{2}{p^{k(P)}}\rf. 
\]

\providecommand{\bysame}{\leavevmode\hbox to3em{\hrulefill}\thinspace}
\providecommand{\MR}{\relax\ifhmode\unskip\space\fi MR }
\providecommand{\MRhref}[2]{%
  \href{http://www.ams.org/mathscinet-getitem?mr=#1}{#2}
}
\providecommand{\href}[2]{#2}

\end{document}